\documentclass[12pt,reqno]{amsart}
\usepackage{hyperref}
\hypersetup{colorlinks,%
	citecolor=black,%
	filecolor=black,%
	linkcolor=black,%
	urlcolor=black,%
	pdftex}
\usepackage{amsthm,amsfonts,amsmath,amssymb,enumerate,longtable}
\newtheorem{thm}{Theorem}

\newtheorem{lem}[thm]{Lemma}

\theoremstyle{definition}

\newtheorem{exam}[thm]{Example}

\numberwithin{equation}{section}
\usepackage{graphicx}
\usepackage[
 a4paper,
 left=2.5cm,
 right=2.5cm,
 top=3cm,
 bottom=3 cm,
 ]{geometry}

\begin{document}
\title{Further comments on a fixed point theorem in the fractal space}

 \author{Nguyen Van Dung}
 
 \address[Nguyen Van Dung]{Faculty of Mathematics Teacher Education, Dong Thap University, Cao Lanh City, Dong Thap Province, Vietnam}
 
 \email{nvdung@dthu.edu.vn}
 
 \author{Wutiphol Sintunavarat}
 
 \address[Wutiphol Sintunavarat]{Department of Mathematics and Statistics, Faculty of Science and Technology, Thammasat University Rangsit Center, 12121 Pathumthani, Thailand}
 
 \email{wutiphol@mathstat.sci.tu.ac.th}

 \thanks{The first author would like to thank the Bualuang ASEAN Fellowship Program for financial support during the preparation of this manuscript.}


\begin{abstract} In this paper, we give some further comments to the counterexample and the results of R.~K. Bisht in [R.~K. Bisht.
	\newblock {Comment on: A new fixed point theorem in the
		fractal space}.
	\newblock {\em Indag. Math. (N.S.)}, 29(2):819--823, 2018.].
\end{abstract}

\maketitle

In \cite{sIR2016}, S. Ri  presented some generalizations of the Banach contraction principle in which the Lipschitz constant $k$ is replaced by some real-valued control function. For the applications to the fractal, the author obtained the fixed point theorem of the some generalized contraction in the fractal space. Then the results have been improved and applied in metric fixed point theory and fractals~\cite{AWS2018}, \cite{JG2017-2}, \cite{KKM2019}, \cite{KS2018-3}, \cite{PP2019}, \cite{tS2016-2}. 
Recently, in \cite{RKB2018-2} R. K. Bisht  gave a counterexample in without proof to \cite[Lemma 2.2]{sIR2016} and improved the result of S. Ri by employing a proper setting. 

In this note, by recalculating the counterexample of R. K. Bisht, we show that the counterexample is not correct. We also present some comments to the main results in~\cite{sIR2016} and~\cite{RKB2018-2}.

For being convenient, we use  the same terminology and the notations as have been utilized in~\cite{sIR2016}. To prove the main result of~\cite{sIR2016}, see \cite[Theorem~2.1]{sIR2016}, S. Ri utilized the following lemma.

\begin{lem}[\cite{sIR2016}, Lemma~2.2] \label{314-00} Assume that the following conditions hold.
\begin{enumerate}
\item $(X,d)$ is a metric space.

\item  $\varphi : [0,\infty) \to [0, \infty)$ is a function with $\varphi(0) =0$, $\varphi (t) <t$ for all $t >0$, $\limsup\limits_{s \to t^+} \varphi (s) < t$ for all $t >0$.
	
\item $f: X \to X$ is a map such that for all $x,y \in X$, 
\begin{equation} \label{314-01} d(f(x),f(y)) \le \varphi (d(x,y)). \end{equation} 
\end{enumerate} Then for each $x \in X$, the sequence $\{f^n(x)\}$ is a Cauchy sequence.
\end{lem}

\begin{thm}[\cite{sIR2016}, Theorem~2.1] \label{314-04} Assume that the following conditions hold.
\begin{enumerate} 
\item \label{314-04-1} $(X,d)$ is a complete metric space.
	
\item \label{314-04-2} $\varphi : [0,\infty) \to [0, \infty)$ is a function with $\varphi(0) =0$, $\varphi (t) <t$ for all $t >0$, $\limsup\limits_{s \to t^+} \varphi (s) < t$ for all $t >0$.

\item $f: X \to X$ is a map such that for all $x,y \in X$, 
\begin{equation} \label{314-03} d(f(x),f(y)) \le \varphi (d(x,y)). \end{equation} 
\end{enumerate} 
Then $f$ has a unique fixed point.
\end{thm}

To disprove Lemma ~\ref{314-00} above, R. K. Bisht gave the following counterexample without proof.

\begin{exam}[\cite{RKB2018-2}, Example~1.2] \label{314-02} 
	Let $X = \Big\{\sum\limits_{ k=1}^n \frac{1}{k}: n=1,2,3,\dots \Big\}$ and $d$ be the usual metric on $X$. Define $f:X \to X$ and $\varphi : [0,\infty ) \to [0,\infty )$ by $\varphi (t) = \frac{t}{1+t}$ for all $t \in [0,\infty )$ and $$f\Big(\sum\limits_{ k=1}^n\frac{1}{k}\Big) = \sum\limits_{ k=1}^{n+1} \frac{1}{k}.$$
Then we have
\begin{enumerate} \item $f$ and $\varphi$ satisfy all the conditions of Lemma~\ref{314-00}.

\item The sequence $\{f^n(x)\}$ is not a Cauchy sequence with $x=1$. 
\end{enumerate} 
\end{exam} 

Unfortunately, we find that for $x = 1$ and $y = 1  + \frac{1}{2} + \frac{1}{3}$,
$$d(f(x),f(y)) = d \left(1  +\frac{1}{2}, 1  + \frac{1}{2} + \frac{1}{3} + \frac{1}{4}  \right) = \frac{1}{3} + \frac{1}{4} = \frac{7}{12}, $$
$$\varphi \left(d(x,y)\right) = \varphi\left(d\left(1, 1  + \frac{1}{2} + \frac{1}{3}\right)\right) = \varphi\left(\frac{1}{2} + \frac{1}{3}\right) = \varphi\left(\frac{5}{6}\right) = \frac{\frac{5}{6}}{1+ \frac{5}{6}} = \frac{5}{11}. $$
Then $d(f(x),f(y)) = \frac{7}{12} > \frac{5}{11} = \varphi(d(x,y)).$ This proves that the condition ~\eqref{314-01} of Lemma~\ref{314-00} does not hold.
Then Example~\ref{314-02} is not correct. However, Theorem~\ref{314-04} still holds. 

R. K. Bisht also proved the following results.

\begin{thm}[\cite{RKB2018-2}, Theorem~2.1]  \label{314-06} Assume that the following conditions hold.
	\begin{enumerate} \item \label{314-06-1} $(X,d)$ is an $f$-orbitally complete metric space.
		
		\item \label{314-06-2} There exists $x_0 \in X$ and  $\varphi_{x_0} : (0,\infty) \to (0, \infty)$ is a function with $\varphi_{x_0} (t) <t$ and $\limsup\limits_{s \to t^+} \varphi_{x_0} (s) < t$ for all $t >0$.
		
		\item $f: X \to X$ is a map such that for all $x, y \in \overline{ O(x_0,f)}$ with $x\neq y$,  
		\begin{eqnarray} \label{314-05} &&d(f(x),f(y)) \\
		& \le&  \varphi_{x_0} \big(\max\{d(x,y), ad(x,f(x)) +(1-a)d(y,f(y)), (1-a)d(x,f(x)) + ad(y,f(y)) \}\big), \nonumber \end{eqnarray} 
		where $\overline{O(x_0,f)}$ is the closure of $O(x_0,f)=\{x_0,fx_0,f^2x_0,f^3x_0,\ldots\}$ and $0<a<1$.
	\end{enumerate} 
	Then we have the following assertions.
\begin{enumerate} 
	\item The sequence $\{f^n(x_0)\}$ is a Cauchy sequence in $X$ and $\lim\limits_{n\to \infty} f^n(x_0) = z \in X$.

\item If $f$ is orbitally continuous at $z$ then $z$ is a  fixed point of $f$.

\item $z$ is the unique fixed point of $f$ in $\overline{ O(x_0,f)}$. \end{enumerate} 
\end{thm} 

\begin{thm}[\cite{RKB2018-2}, Theorem~2.3]  \label{314-09} Assume that the following conditions hold. 
	\begin{enumerate} \item \label{314-09-1} $(X,d)$ is an $f$-orbitally complete metric space.
		
		\item \label{314-09-2} There exists $x_0 \in X$ and  $\varphi_{x_0} : (0,\infty) \to (0, \infty)$ is a function with $\varphi_{x_0} (t) <t$ and $\limsup\limits_{s \to t^+} \varphi_{x_0} (s) < t$ for all $t >0$.
		
		\item $f: X \to X$ is a map such that for all $x, y \in \overline{ O(x_0,f)}$ with $x\neq y$, 
		\begin{eqnarray} \label{314-10} d(f(x),f(y))
		& \le&  \varphi_{x_0} \big(\max\{d(x,y),d(x,f(x)),d(y,f(y)) \}\big),  \end{eqnarray} 
	\end{enumerate} 
	where $\overline{O(x_0,f)}$ is the closure of $O(x_0,f)=\{x_0,fx_0,f^2x_0,f^3x_0,\ldots\}$.
	Then we have the following assertions.
	\begin{enumerate} 
		\item The sequence $\{f^n(x_0)\}$ is a Cauchy sequence  in $X$ and $\lim\limits_{n\to \infty} f^n(x_0) = z \in X$.
		
		\item If $f$ is orbitally continuous at $z$ then $z$ is a  fixed point of $f$.
		
		\item $z$ is the unique fixed point of $f$ in $\overline{ O(x_0,f)}$. \end{enumerate} 
\end{thm} 

We have some comments on Theorem~\ref{314-06}  and Theorem~\ref{314-09} as follows. 

\begin{enumerate} \item Theorem~\ref{314-04} assumes the conditions for the complete metric space $X$ while Theorem~\ref{314-06} and Theorem~\ref{314-09} assume the condition for the complete metric space $\overline{O(x_0,f)}$. The calculations are the same. This idea first appeared in ~\cite{LBC1974}.
	
\item The function $\varphi $ is from $[0,\infty )$ to $[0,\infty )$ in Theorem~\ref{314-04} and the function $\varphi_{x_0}$ is from $(0,\infty )$ to $(0,\infty )$ in Theorem~\ref{314-06} and Theorem~\ref{314-09}, then the condition~\eqref{314-05} and the condition~\eqref{314-10} are for $x \ne y$ to satisfy that $\varphi_{x_0}$ is not defined at 0.

\item The assumption of orbital continuity at $z$ of $f$ in Theorem~\ref{314-06} is redundant. Indeed, from~\eqref{314-05} and $\varphi_{x_0}(t) <t$ for all $t >0$, we have for all $x, y \in \overline{ O(x_0,f)}$, 
\begin{equation*}  d(f(x),f(y)) \le \max\{d(x,y), ad(x,f(x)) +(1-a)d(y,f(y)), (1-a)d(x,f(x)) + ad(y,f(y)) \}. \end{equation*}  Note that $z \in \overline{ O(x_0,f)}$. So we have
\begin{eqnarray} \label{314-07} d(f^{n+1}(x_0),f(z)) 
& \le &\max\{d(f^n(x_0),z), ad(f^n(x_0),f^{n+1}(x_0)) +(1-a)d(z,f(z)),\nonumber\\
&& (1-a)d(f^n(x_0),f^{n+1}(x_0)) + ad(z,f(z)) \}. \end{eqnarray}
Letting $n \to  \infty $ in ~\eqref{314-07} and using $\lim\limits_{n\to \infty} f^n(x_0) = z$ we have 
\begin{equation} \label{314-08} d(z,f(z)) \le \max\{(1-a)d(z,f(z)), ad(z,d(z))\} = \max\{(1-a), a\} d(z,f(z)) . \end{equation}  Note that $0<a<1$. So from ~\eqref{314-08}  we get $d(z,f(z)) = 0$, that is, $z$ is a  fixed point of $f$.

\item The assumption of orbital continuity at $z$ of $f$ in Theorem~\ref{314-09} is also redundant. Indeed, if there exists $n_0$ such that $f^n(x_0) = z$ for all $n \ge n_0$, then $z$ is a fixed point of $f$. Otherwise, the exists a subsequence $\{f^{k_n}(x_0)\}$ of $\{f^n(x_0)\}$ such that $f^{k_n}(x_0) \ne z$ for all $k_n$. Moreover, the subsequence can be chosen such that the sequence $\{d(f^{k_n}(x_0),z)\}$ is decreasing to 0.  Then, from~\eqref{314-10} and $z \in \overline{ O(x_0,f)}$, we have
\begin{eqnarray} \label{314-12} d(f^{k_n+1}(x_0),f(z)) 
& \le &\varphi_{x_0}\big(\max\{d(f^{k_n}(x_0),z), d(f^{k_n}(x_0),f^{k_n+1}(x_0)),d(z,f(z)) \}\big). \end{eqnarray}
Note that the sequence $\big\{ \max\{d(f^{k_n}(x_0),z), d(f^{k_n}(x_0),f^{k_n+1}(x_0)),d(z,f(z)) \}\big\}$ is decreasing to $d(z,f(z))$. Suppose to the contrary that $d(z,f(z)) > 0$. Then letting $n \to  \infty $ in ~\eqref{314-12} and using $\limsup\limits_{s \to t^+} \varphi_{x_0} (s) < t$ for all $t >0$, we have 
\begin{eqnarray*} \label{314-13} d(z,f(z)) &\le & \limsup\limits_{ n \to \infty }\varphi_{x_0}\big(\max\{d(f^{k_n}(x_0),z), d(f^{k_n}(x_0),f^{k_n+1}(x_0)),d(z,f(z)) \}\big) \\
	& <&  d(z,f(z)). \end{eqnarray*} This is a contradiction. So $d(z,f(z)) = 0$, $z$ is a  fixed point of $f$. 
 \end{enumerate}

\section*{Acknowledgements}

The second author would like to thank the Thailand Research Fund and Office of the Higher Education Commission under grant no. MRG6180283 for financial support during the preparation of this manuscript.

\bibliographystyle{abbrv}  

\end{document}